\documentclass[11pt]{article}
\usepackage{amsmath,amssymb,latexsym,epsf,mathabx,color}
\usepackage[all]{xy}

\begin{document}

\newtheorem{Th}{Theorem}[section]
\newtheorem{Def}[Th]{Definition}
\newtheorem{Lem}[Th]{Lemma}
\newtheorem{Pro}[Th]{Proposition}
\newtheorem{Cor}[Th]{Corollary}
\newtheorem{Rem}[Th]{Remark}
\newtheorem{Exm}[Th]{Example}
\newtheorem{Sc}[Th]{}
\def\Pf#1{{\noindent\bf Proof}.\setcounter{equation}{0}}
\def\bg#1{\begin{#1}\setcounter{equation}{0}}
\def\ed#1{\end{#1}\setcounter{equation}{0}}


\title{\bf  Silting objects and torsion pairs in comma categories\thanks{Supported by the Natural Science Foundation of Universities of Anhui (No.2023AH050950 and 2023AH050904), the Top talent project of AHPU in 2020 (No.S022021055),
the National Natural Science Foundation of China (Grants No.11801004, 12101003 and 12301042), the Natural Science Foundation of Anhui province (No.2108085QA07) and
the Startup Foundation for Introducing Talent of AHPU (No.2020YQQ067 and 2022YQQ09).}}

\smallskip
\author{ Peiyu Zhang\footnote{Correspondent author}, Xinyu Wang, Dajun Liu and Li Wang\\ 
\footnotesize ~E-mail:~zhangpy@ahpu.edu.cn,  w13955807808@163.com, liudajun@ahpu.edu.cn and wangli1995@ahpu.edu.cn; \\
\footnotesize School of Mathematics-Physics and Finance,  Anhui Polytechnic University, Wuhu, China.\\
Jiaqun Wei\\
\footnotesize ~E-mail: weijiaqun@njnu.edu.cn \\
\footnotesize  School of Mathematical Sciences, Zhejiang Normal University, Jinhua, China\\
}
\date{}
\maketitle
%
%
%
\begin{abstract}
\vskip 10pt%
In this paper, firstly, we give a characterization of silting objects in the comma category $(F,~\mathrm{Mod} S)$.
Let $U$ be an $(S,~R)$-bimodule. Take $F=U\bigotimes_{R}-$, then the comma category $(F,~\mathrm{Mod} S)$ is isomorphic to the category $T$-Mod, where
$T$ is a triangular matrix ring. Assume that $\mathcal{C}_{1}$ and $\mathcal{C}_{2}$ are two subcategories of left $R$-modules,
$\mathcal{D}_{1}$ and $\mathcal{D}_{2}$ are two subcategories of left $S$-modules.
We prove that $(\mathcal{C}_{1},~\mathcal{C}_{2})$ and $(\mathcal{D}_{1},~\mathcal{D}_{2})$ are torsion pairs if and only if
$(\mathfrak{B}^{\mathcal{C}_{1}}_{\mathcal{D}_{1}},~\mathfrak{U}^{\mathcal{C}_{2}}_{\mathcal{D}_{2}})$ and
$(\mathfrak{U}^{\mathcal{C}_{1}}_{\mathcal{D}_{1}},~\mathfrak{J}^{\mathcal{C}_{2}}_{\mathcal{D}_{2}})$ are torsion pairs under certain conditions.

\vskip 10pt

\noindent MR(2020) Subject Classification: 16D90; 16E30; 16S90.



\end{abstract}
%
\vskip 30pt

\section{Introduction and Preliminaries}

%
%
%
\hskip 18pt
The tilting theory is well known, and plays an important role in the representation theory of
Artin algebra. The classical notion of tilting modules was considered in the case
of finite dimensional algebras by Brenner and Butler \cite{SBMB} and by Happel and Ringel \cite{HR}.
As a generalization of tilting modules, silting modules over arbitrary rings  were introduced by
L. Angeleri-H\"{u}gel, F. Marks and J. Vit\'{o}ria, and they proved that silting modules are in bijection with
2-term silting complexes and with certain t-structures and co-t-structures in the derived module category.

Let $\mathcal{A}$ be an abelian category. Given a subcategory $\mathcal{X}$ of $\mathcal{A}$, write $\mathcal{X}^{\perp_{0}}=\{Y: \mathrm{Hom}_{\mathcal{A}}(X,~Y)=0$ for all $X\in \mathcal{X}\}$
and $^{\perp_{0}}\mathcal{X}=\{Y: \mathrm{Hom}_{\mathcal{A}}(Y,~X)=0$ for all $X\in \mathcal{X}\}$.
Recall that a pair of subcategories $(\mathcal{X}, ~\mathcal{Y })$ is called a torsion pair \cite{ASS,D} if the following conditions are satisfied:
(1) $\mathrm{Hom}_{\mathcal{A}}(\mathcal{X},~\mathcal{Y})= 0$;
(2) for any object $M\in\mathcal{A}$, there exists an exact sequence $\xymatrix{0\ar[r]&X\ar[r]&M\ar[r]&Y\ar[r]&0}$
in $\mathcal{A}$ with $X\in\mathcal{X}$ and $Y\in\mathcal{Y}$.
It is equivalent to that $\mathcal{X}^{\perp_{0}}=\mathcal{Y}$ and $^{\perp_{0}}\mathcal{Y}=\mathcal{X}$.
As we all known, for any silting (resp. tilting) module $T$, it can induce a torsion pair $(\mathrm{Gen}T,~T^{\perp_{0}})$.

Recall that for two Abelian categories $\mathcal{A}$ and $\mathcal{B}$, and a right exact functor
$F$: $\mathcal{A}\rightarrow \mathcal{B}$, we can define the left comma category $(F,~\mathcal{B})$ \cite{FGR},
which is also an Abelian category. The example of comma categories include but
are not limited to the category of modules or complexes over a triangular matrix ring, the morphism
category of an abelian category and so on \cite{FGR,MLX,HZ}.
Other studies of comma categories, readers can also refer to \cite{CL,L,MAR}.

This paper is organized as follows. In Section 2, we firstly introduce some definitions and properties,
and give a characterization of silting objects in the comma category $(F,~\mathrm{Mod} S)$, see Theorem \ref{1}.
In section 3, we study three functors $(\textbf{p},~ \textbf{q},~\textbf{h})$ and three classes,
and mainly prove Theorem \ref{2}.

\bg{Th}\rm{}\label{1}
Let $A\in R$-$\mathrm{Mod}$, $B\in S$-$\mathrm{Mod}$ and the functor $F$ commute with direct sums.
Then $\binom{A}{FA\bigoplus B}_{\binom{1}{0}}$ is silting with respect to $\sigma$ if and only if

$(1)$ $A$ is silting with respect to $\sigma_{A}$;

$(2)$ $B$ is silting with respect to $\sigma_{B}$;

$(3)$ $FA\in \mathrm{Gen}B$.
\ed{Th}

\bg{Th}\rm{}\label{2}
Let $\mathcal{C}_{1}$ and $\mathcal{C}_{2}$ be two subcategories of left $R$-modules, $\mathcal{D}_{1}$ and $\mathcal{D}_{2}$ be two subcategories of left $S$-modules,
$_{R}U$ be flat and $U_{S}$ be projective.

(1) If $S^{+}\in\mathcal{D}_{2}$, then $(\mathcal{C}_{1},~\mathcal{C}_{2})$ and $(\mathcal{D}_{1},~\mathcal{D}_{2})$ are torsion pairs if and only if
$(\mathfrak{B}^{\mathcal{C}_{1}}_{\mathcal{D}_{1}},~\mathfrak{U}^{\mathcal{C}_{2}}_{\mathcal{D}_{2}})$ is a torsion pair.

(2) If $R\in\mathcal{C}_{1}$, then $(\mathcal{C}_{1},~\mathcal{C}_{2})$ and $(\mathcal{D}_{1},~\mathcal{D}_{2})$ are torsion pairs if and only if
$(\mathfrak{U}^{\mathcal{C}_{1}}_{\mathcal{D}_{1}},~\mathfrak{J}^{\mathcal{C}_{2}}_{\mathcal{D}_{2}})$ is a torsion pair.
\ed{Th}



\section{Silting objects}
In this section, we firstly recall some definitions and notions. Let $\mathcal{A}$ be an abelian category.
Recall that a subcategory $\mathcal{C}\subseteq\mathcal{A}$ is called a torsion class, if it is closed under images, direct sums and extensions (c.f. \cite[Chapter~VI]{ASS}).
Given a subcategory $\mathcal{X}\subseteq\mathcal{A}$, recalled that a left $\mathcal{X}$-approximation of $T\in\mathcal{A}$ is a morphism $\phi:$ $T\longrightarrow X$
such that $\mathrm{Hom}_{\mathcal{A}}(\phi,~X')$ is surjective for any $X'\in \mathcal{X}$.

For a morphism $f$ in $\mathrm{Proj}\mathcal{A}$ consisting of all projective objects of $\mathcal{A}$, we consider the following class of objects
$$\mathcal{D}_{f}=\{X\in \mathcal{A}|\mathrm{Hom}_{\mathcal{A}}(f,X)\rm{~is~ surjective}\}.$$

\bg{Def}\rm{\cite[Definition~3.7]{AMV}}\label{silting}%
Let $\mathcal{A}$ be an abelian category with enough projective objects and $T$ be an object of $\mathcal{A}$.

$(1)$ $T$ is called partial silting if there exists a projective presentation $\sigma$ of $T$ such that $\mathcal{D}_{\sigma}$ is a torsion class and $T\in\mathcal{D}_{\sigma}$.

$(2)$ $T$ is called silting if there exists a projective presentation $\sigma$ of $T$ such that $\mathcal{D}_{\sigma}=\mathrm{Gen}T$,
where $\mathrm{Gen}T=\{X|$ there is a surjective morphism $T^{(I)}\longrightarrow X$, for some set $I$ $\}$.

In general, we also say that $T$ is (partial) silting with respect to $\sigma$.
\ed{Def}

By \cite[Lemma~3.6]{AMV}, $\mathcal{D}_{\sigma}$ is always closed under images and extensions. Hence $\mathcal{D}_{\sigma}$ is a torsion class if and only if it is colsed under direct sums.

\bg{Def}\rm(\cite[Section~1]{FGR})
Let $\mathcal{A}$ be an abelian category and $F$: $\mathcal{A}\rightarrow \mathcal{A}$ be an additive endofunctor. The right
trivial extension of $\mathcal{A}$ by $F$, denoted by $\mathcal{A}\ltimes F$, is defined as follows. An object in $\mathcal{A}\ltimes F$ is a morphism
$\alpha$: $FA\rightarrow A$ for an object $A$ in $\mathcal{A}$ such that $\alpha\cdot F(\alpha)=0$; and a morphism in $\mathcal{A}\ltimes F$ is a pair $(F\beta,~\beta)$
of morphisms in $\mathcal{ A}$ such that the following diagram
$$\xymatrix{
FA\ar[d]^{F\beta}\ar[r]^{\alpha}&A\ar[d]^{\beta}\\
FA'\ar[r]^{\alpha'}&A'
}$$
is commutative.
\ed{Def}

\bg{Def}\rm(\cite[Section~1]{FGR})\label{comma}
Let $\mathcal{A}$ and $\mathcal{B}$ be two abelian categories and $F$: $\mathcal{A}\rightarrow \mathcal{B}$ be an additive functor.
We define the left comma category $(F,~\mathcal{B})$ as follows. The objects of the category are $\binom{A}{B}_{\phi}$ (sometimes, $\phi$ is omitted) with
$A\in \mathcal{A}$, $B\in \mathcal{B}$ and $\phi\in \mathrm{Hom}(FA,~B)$;
and the morphisms are pairs $\binom{\alpha}{\beta}$: $\binom{A}{B}_{\phi}\longrightarrow\binom{A'}{B'}_{\phi'}$  with $\alpha\in \mathrm{Hom}_{\mathcal{A}}(A,~A')$ and $\beta\in \mathrm{Hom}_{\mathcal{B}}(B,~B')$
such that that the following diagram
$$\xymatrix{
FA\ar[d]^{F\alpha}\ar[r]^{\phi}&B\ar[d]^{\beta}\\
FA'\ar[r]^{\phi'}&B'
}$$
is commutative.
\ed{Def}

Similarly, we have the right comma category $(G,~\mathcal{A})$, where $G$: $\mathcal{B}\longrightarrow \mathcal{A}$ is an additive functor.
Its objects are triples $(X,~Y)_{\phi}$, where $X\in\mathcal{A}$, $Y\in\mathcal{B}$ and $\phi$: $GY\longrightarrow X$,
and whose morphisms from $(X_{1},~Y_{1})_{\phi_{1}}$ to $(X_{2},~Y_{2})_{\phi_{2}}$ are pairs $(f,~g)$ such that $f\in\mathrm{Hom}_{\mathcal{A}}(X_{1},~X_{2})$, $g\in\mathrm{Hom}_{\mathcal{B}}(Y_{1},~Y_{2})$
and $\phi_{2}Gg=f\phi_{1}$. In this paper, we mainly consider the left comma category.

Recall that the product category $\mathcal{A}\times\mathcal{B}$ is defined as follows: an object of $\mathcal{A}\times\mathcal{B}$
is a pair $(A,~B)$ with $A\in \mathcal{A}$ and $B\in \mathcal{B}$, a morphism from $(A,~B)$ to $(A',~B')$ is a pair $(f,~g)$ with
$f\in\mathrm{Hom}_{\mathcal{A}}(A,~A')$ and $g\in\mathrm{Hom}_{\mathcal{B}}(B,~B')$.

Next, we give some comments on the comma category, which are very helpful for us to understand the comma category.

\bg{Rem}\rm{}\label{rem1}
Let $\mathcal{A}$ and $\mathcal{B}$ be two abelian categories and $F$: $\mathcal{A}\rightarrow \mathcal{B}$ be an additive functor.

$(1)$ The functor $F$ induces a functor $G$: $\mathcal{A}\times\mathcal{B}\longrightarrow\mathcal{A}\times\mathcal{B}$ by $G(A,B)=(0,FA)$ and
$G(\alpha,\beta)=(0,F\alpha)$. It is easy to verify that $(F,\mathcal{B})\cong (\mathcal{A}\times\mathcal{B})\ltimes G$, by mapping the object $\binom{A}{B}_{\phi}$
to the object $(0,\phi)$: $G(A,B)\longrightarrow(A,B)$ in $(\mathcal{A}\times\mathcal{B})\ltimes G$.

$(2)$ Suppose that $F$ is right exact. It is easy to see that the above functor $G$ is also right exact. So $(\mathcal{A}\times\mathcal{B})\ltimes G$ is an abelian
category by \cite[Proposition~1.1]{FGR}, and then the comma category $(F,\mathcal{B})$ is also abelian by (1).

$(3)$ The sequence $\xymatrix{0\ar[r]&\binom{A_{1}}{B_{1}}_{\phi_{1}}\ar[r]&\binom{A_{2}}{B_{2}}_{\phi_{2}}\ar[r]&\binom{A_{3}}{B_{3}}_{\phi_{3}}\ar[r]&0}$ is exact in $(F,\mathcal{B})$
if and only if the sequence $\xymatrix{0\ar[r]&(0,~\phi_{1})\ar[r]&(0,~\phi_{2})\ar[r]&(0,~\phi_{3})\ar[r]&0}$ is exact in $(\mathcal{A}\times\mathcal{B})\ltimes G$ by (1)
if and only if the sequence $\xymatrix{0\ar[r]&(A_{1},~B_{1})\ar[r]&(A_{2},~B_{2})\ar[r]&(A_{3},~B_{3})\ar[r]&0}$ is exact in $\mathcal{A}\times\mathcal{B}$ by \cite[Corollary~1.2]{FGR}
if and only if the following two sequences $\xymatrix@C=0.5cm{0\ar[r]&A_{1}\ar[r]&A_{2}\ar[r]&A_{3}\ar[r]&0}$ in $\mathcal{A}$ and
$\xymatrix{0\ar[r]&B_{1}\ar[r]&B_{2}\ar[r]&B_{3}\ar[r]&0}$ in $\mathcal{B}$ are exact.

$(4)$ The projective object in $(F,\mathcal{B})$ is of the form $\binom{0}{Q}\bigoplus\binom{P}{FP}$ with $P$ projective in $\mathcal{A}$ and $Q$ projective in $\mathcal{B}$ by \cite[Lemma~3.1]{PZH}.
\ed{Rem}

From now on, $R$ and $S$ are rings, denoted $R$-$\mathrm{Mod}$ by the category consisting of all left $R$-modules. The functor $F$: $R$-$\mathrm{Mod}\longrightarrow$ $S$-$\mathrm{Mod}$ is right exact and covariant.
In this subsection, we mainly consider the left comma category $\mathcal{C}:=(F$, $S$-$\mathrm{Mod})$.

Let $A\in R$-$\mathrm{Mod}$ and $B\in S$-$\mathrm{Mod}$, and the following sequences
$$\xymatrix{P_{1}\ar[r]^{\sigma_{A}}&P_{0}\ar[r]&A\ar[r]&0}$$
and
$$\xymatrix{Q_{1}\ar[r]^{\sigma_{B}}&Q_{0}\ar[r]&B\ar[r]&0}$$
are projective presentations of $A$ and $B$, respectively, with $P_{0},~P_{1}$ being projective $R$-modules and $Q_{0},~Q_{1}$ being projective $S$-modules.
Since $F$ is right exact and covariant, the sequence $\xymatrix{FP_{1}\ar[r]^{F\sigma_{A}}&FP_{0}\ar[r]&FA\ar[r]&0}$ is also exact.
Hence, by Remark \ref{rem1}(4), we can obtain a projective presentation of $\binom{0}{B}\bigoplus\binom{A}{FA}$, denoted by $\sigma=\begin{pmatrix}a&0\\0&b\end{pmatrix}$:
$$\xymatrix{\binom{0}{Q_{1}}\bigoplus\binom{P_{1}}{FP_{1}}\ar[r]^{\sigma}&\binom{0}{Q_{0}}\bigoplus\binom{P_{0}}{FP_{0}}\ar[r]&\binom{0}{B}\bigoplus\binom{A}{FA}\ar[r]&0,}$$
where $a=\binom{0}{\sigma_{B}}$, $b=\binom{\sigma_{A}}{F\sigma_{A}}$.

\vskip 10pt
Next, we will give the characterization of $\mathcal{D}_{\sigma}$ in the comma category.

\bg{Lem}\label{lem1}%
$\binom{M}{N}_{h}\in \mathcal{D}_{\sigma}$ if and only if $M\in \mathcal{D}_{\sigma_{A}}$ and $N\in \mathcal{D}_{\sigma_{B}}$.
\ed{Lem}

\Pf. ($\Rightarrow$) For any morphism $f\in \mathrm{Hom}_{R}(P_{1},~M)$ and $m\in\mathrm{Hom}_{S}(Q_{1},~N)$, we only need to find two morphisms such that the following diagrams are commutative.
$$\xymatrix{
M&\\
P_{1}\ar[r]^{\sigma_{A}}\ar[u]^{f}&P_{0}\ar@{.>}[ul]_{?}\ar[r]&0
}~~~~~~
\xymatrix{
N&\\
Q_{1}\ar[r]^{\sigma_{B}}\ar[u]^{m}&Q_{0}\ar@{.>}[ul]_{?}\ar[r]&0
}$$
Set $g=h\cdot Ff$, then the following diagram is commutative.
$$\xymatrix{
FP_{1}\ar[d]^{Ff}\ar[r]^{1}&FP_{1}\ar[d]^{g}\\
FM\ar[r]^{h}&N
}$$
i.e., $(\binom{0}{m},~\binom{f}{g})\in \mathrm{Hom}_{\mathcal{C}}(\binom{0}{Q_{1}}\bigoplus\binom{P_{1}}{FP_{1}},~\binom{M}{N})$.
By the assumption, there is a morphism $(\binom{0}{i},~\binom{x}{y})$: $\binom{0}{Q_{0}}\bigoplus\binom{P_{0}}{FP_{0}}\longrightarrow\binom{M}{N}$
such that the following diagram is commutative.
$$\xymatrix{
\binom{M}{N}_{h}&\\
\binom{0}{Q_{1}}\bigoplus\binom{P_{1}}{FP_{1}}\ar[u]^{(\binom{0}{m},~\binom{f}{g})}\ar[r]^{\sigma}&\binom{0}{Q_{0}}\bigoplus\binom{P_{0}}{FP_{0}}\ar@{.>}[ul]_{(\binom{0}{i},~\binom{x}{y})}
}$$
By the definition of morphisms in the comma category, it is easy to prove that $y=h\cdot Fx$. Since the above diagram is commutative, we have the following two equations
$\binom{0}{m}=\binom{0}{i}\binom{0}{\sigma_{B}}$ and $\binom{f}{g}=\binom{x}{y}\binom{\sigma_{A}}{F\sigma_{A}}$.
Thus $f=x\sigma_{A}$ and $m=i\sigma_{B}$. i.e., $M\in \mathcal{D}_{\sigma_{A}}$ and $N\in \mathcal{D}_{\sigma_{B}}$.

($\Leftarrow$) For any $(\binom{0}{j},~\binom{c}{d})\in \mathrm{Hom}_{\mathcal{C}}(\binom{0}{Q_{1}}\bigoplus\binom{P_{1}}{FP_{1}},~\binom{M}{N}_{h})$.
By the definition of morphisms in the comma category, it is easy to verify that $d=h\cdot Fc$. Now, we need to find a morphism such that the following diagram is commutative.
$$\xymatrix{
\binom{M}{N}_{h}&\\
\binom{0}{Q_{1}}\bigoplus\binom{P_{1}}{FP_{1}}\ar[u]^{(\binom{0}{j},~\binom{c}{d})}\ar[r]^{\sigma}&\binom{0}{Q_{0}}\bigoplus\binom{P_{0}}{FP_{0}}\ar@{.>}[ul]_{?}
}$$
Since $M\in \mathcal{D}_{\sigma_{A}}$ and $N\in \mathcal{D}_{\sigma_{B}}$, there are two morphisms $j_{1}$ and $c_{1}$ satisfying the following two diagrams are commutative.
$$\xymatrix{
M&\\
P_{1}\ar[r]^{\sigma_{A}}\ar[u]^{c}&P_{0}\ar@{.>}[ul]_{\exists ~c_{1}}\ar[r]&0
}~~~~~~
\xymatrix{
N&\\
Q_{1}\ar[r]^{\sigma_{B}}\ar[u]^{j}&Q_{0}\ar@{.>}[ul]_{\exists ~j_{1}}\ar[r]&0
}$$
Set $d_{1}=h\cdot Fc_{1}$, then we have the following commutative diagram.
$$\xymatrix{
FP_{0}\ar[d]^{1}\ar[r]^{Fc_{1}}&FM\ar[d]^{h}\\
FP_{0}\ar[r]^{d_{1}}&N
}$$
We claim that the morphism $(\binom{0}{j_{1}},~\binom{c_{1}}{d_{1}})\in \mathrm{Hom}_{\mathcal{C}}(\binom{0}{Q_{0}}\bigoplus\binom{P_{0}}{FP_{0}},~\binom{M}{N}_{h})$ is exactly what we want.
Thus we need to verify that the following two equations hold.
$$\binom{0}{j}=\binom{0}{j_{1}}\binom{0}{\sigma_{B}}~~~~\mathrm{and}~~~~
\binom{c}{d}=\binom{c_{1}}{d_{1}}\binom{\sigma_{A}}{F\sigma_{A}}.$$
Note that $d=h\cdot Fc=h\cdot F(c_{1}\sigma_{A})=(h\cdot F(c_{1}))\cdot F(\sigma_{A})=d_{1}\cdot F(\sigma_{A})$. Thus the second equation holds.
The first is obvious. So we completed this proof.
\ \hfill $\Box$

\vskip 10pt
By the above conclusion, we obtain the following corollary immediately.

\bg{Cor}\rm{}\label{cor1}%
Let $M\in R$-$\mathrm{Mod}$ and $N\in S$-$\mathrm{Mod}$, then the following statements hold.

$(1)$ If $M\in\mathcal{D}_{\sigma_{A}}$, then $\binom{M}{0}\in \mathcal{D}_{\sigma}$.

$(2)$ If $N\in\mathcal{D}_{\sigma_{B}}$, then $\binom{0}{N}\in \mathcal{D}_{\sigma}$.

$(3)$ If $M\in\mathcal{D}_{\sigma_{A}}$ and $FM\in\mathcal{D}_{\sigma_{B}}$, then $\binom{M}{FM}\in \mathcal{D}_{\sigma}$.
\ed{Cor}

\bg{Lem}\rm{}\label{lem2}%
Let $A\in R$-$\mathrm{Mod}$, $B\in S$-$\mathrm{Mod}$ and the functor $F$ commutes with direct sums.
Then ${D}_{\sigma}$ is a torsion class if and only if both $\mathcal{D}_{\sigma_{A}}$ and $\mathcal{D}_{\sigma_{B}}$ are torsion class.
\ed{Lem}

\Pf. By \cite[Lemma~3.6]{AMV}, $\mathcal{D}_{\sigma}$ is always closed under images and extensions. Hence $\mathcal{D}_{\sigma}$ is a torsion class if and only if it is colsed under direct sums.

($\Rightarrow$) For any $M_{i}\in\mathcal{D}_{\sigma_{A}}$, then $\binom{ M_{i}}{0}\in\mathcal{D}_{\sigma}$ by Corollary \ref{cor1}.
Since ${D}_{\sigma}$ is a torsion class and $F$ commute with direct sums, $\bigoplus_{i\in I}\binom{ M_{i}}{0}\cong \binom{\bigoplus_{i\in I} M_{i}}{0}\in\mathcal{D}_{\sigma}$ for any set $I$.
By Lemma \ref{lem1}, $\bigoplus_{i\in I}M_{i}\in\mathcal{D}_{\sigma_{A}}$ for any set $I$. i.e., $\mathcal{D}_{\sigma_{A}}$ is a torsion class.
Similarly, we can prove that  $\mathcal{D}_{\sigma_{B}}$ is a torsion class.

($\Leftarrow$) For any $\binom{ X_{i}}{Y_{i}}_{f_{i}}\in\mathcal{D}_{\sigma}$, then $X_{i}\in\mathcal{D}_{\sigma_{A}}$ and $Y_{i}\in\mathcal{D}_{\sigma_{B}}$ by Lemma \ref{lem1}.
Since the functor $F$ be commutes with direct sums and $\mathcal{D}_{\sigma}$ is a torsion class, we have that
$\bigoplus\binom{ X_{i}}{Y_{i}}_{f_{i}}\cong\binom{ \bigoplus X_{i}}{\bigoplus Y_{i}}_{\bigoplus f_{i}}$. Since both $\mathcal{D}_{\sigma_{A}}$ and $\mathcal{D}_{\sigma_{B}}$ are torsion class,
$\bigoplus X_{i}\in\mathcal{D}_{\sigma_{A}}$ and $\bigoplus Y_{i}\in\mathcal{D}_{\sigma_{B}}$.
By Lemma \ref{lem1}, $\bigoplus\binom{ X_{i}}{Y_{i}}_{f_{i}}\cong\binom{ \bigoplus X_{i}}{\bigoplus Y_{i}}_{\bigoplus f_{i}}\in\mathcal{D}_{\sigma}$,
i.e., ${D}_{\sigma}$ is a torsion class.
\ \hfill $\Box$

\bg{Pro}\rm{}\label{prop1}%
Let $A\in R$-$\mathrm{Mod}$ and $B\in S$-$\mathrm{Mod}$ and the functor $F$ commute with direct sums.
Then $\binom{A}{FA}\bigoplus\binom{0}{B}$ is partial silting with respect to $\sigma$ if and only if

$(1)$ $A$ is partial silting with respect to $\sigma_{A}$;

$(2)$ $B$ is partial silting with respect to $\sigma_{B}$;

$(3)$ $FA\in \mathcal{D}_{\sigma_{B}}$.
\ed{Pro}

\Pf. By Lemma \ref{lem1}, Corollary \ref{cor1} and Lemma \ref{lem2}, it is obvious.
\ \hfill $\Box$

\bg{Lem}\rm{\cite[Proposition~3.11]{AMV}}\label{lem3}%
Let $T\in R$-$\mathrm{Mod}$ with a projective presentation $\sigma_{T}$. Then $T$ is silting with respect to $\sigma_{T}$
if and only if $T$ is partial silting with respect to $\sigma_{T}$ and there exists an exact sequence
$$\xymatrix{R\ar[r]^{\alpha}&T_{0}\ar[r]&T_{1}\ar[r]&0}$$
in $\mathrm{Mod}R$ with $T_{0}$, $T_{1}\in \mathrm{Add}T$ and $\alpha$ a left $\mathcal{D}_{\sigma_{T}}$-approximation.
\ed{Lem}

\vskip 10pt

Next, we will give the following main result in this section.

\bg{Th}\rm{}\label{th}%
Let $A\in R$-$\mathrm{Mod}$, $B\in S$-$\mathrm{Mod}$ and the functor $F$ commute with direct sums.
Then $\binom{A}{FA\bigoplus B}_{\binom{1}{0}}$ is silting with respect to $\sigma$ if and only if

$(1)$ $A$ is is silting with respect to $\sigma_{A}$;

$(2)$ $B$ is is silting with respect to $\sigma_{B}$;

$(3)$ $FA\in \mathrm{Gen}B$.
\ed{Th}

\Pf. Noth that $\binom{A}{FA\bigoplus B}_{\binom{1}{0}}\cong\binom{A}{FA}\bigoplus\binom{0}{B}$.

($\Rightarrow$) Since $\binom{A}{FA}\bigoplus\binom{0}{B}$ is silting,
$\mathcal{D}_{\sigma}=\mathrm{Gen}(\binom{A}{FA}\bigoplus\binom{0}{B})=\mathrm{Gen}\binom{A}{FA}\bigoplus \mathrm{Gen}\binom{0}{B}$.
By Proposition $\ref{prop1}$, $A$ is partial silting, and then $\mathrm{Gen}A\subseteq \mathcal{D}_{\sigma_{A}}$.
Next, we will prove that $\mathcal{D}_{\sigma_{A}}\subseteq\mathrm{Gen}A$. For any $X\in \mathcal{D}_{\sigma_{A}}$,
by Corollary \ref{cor1}, we have that $\binom{X}{0}\in \mathcal{D}_{\sigma}=\mathrm{Gen}\binom{A}{FA}\bigoplus \mathrm{Gen}\binom{0}{B}$.
We claim that $\binom{X}{0}\in\mathrm{Gen}\binom{A}{FA}$. If $\binom{X}{0}$ has a direct summand $\binom{X_{1}}{0}\in\mathrm{Gen}\binom{0}{B}$, then there is
a exact sequence $\xymatrix{\binom{0}{B^{(I)}}\ar[r]&\binom{X_{1}}{0}\ar[r]&0}$,  i.e., $\xymatrix{0\ar[r]&X_{1}\ar[r]&0}$ is exact by Remark \ref{rem1},
and then $X_{1}=0$. So we have that $\binom{X}{0}\in\mathrm{Gen}\binom{A}{FA}$. By Remark \ref{rem1} again, we have that $X\in\mathrm{Gen}A$,
and then $\mathcal{D}_{\sigma_{A}}\subseteq\mathrm{Gen}A$. Consequently, $A$ is silting with respect to $\sigma_{A}$.
Similarly, we can prove that $B$ is silting with respect to $\sigma_{B}$. Note that $\binom{A}{FA}\in\mathrm{Gen}(\binom{A}{FA}\bigoplus\binom{0}{B})=\mathcal{D}_{\sigma}$.
Then we have that $FA\in\mathcal{D}_{\sigma_{B}}=\mathrm{Gen}B$ by Lemma \ref{lem1}.

($\Leftarrow$) Note that $\binom{A}{FA}\bigoplus\binom{0}{B}$ is partial silting by Proposition \ref{prop1}, and then $\mathrm{Gen}(\binom{A}{FA}\bigoplus\binom{0}{B})\subseteq\mathcal{D}_{\sigma}$
by Remark 3.8 in \cite{AMV}. So we only need to prove that $\mathcal{D}_{\sigma}\subseteq\mathrm{Gen}(\binom{A}{FA}\bigoplus\binom{0}{B})$.

Since both $A$ and $B$ are silting in $R$-$\mathrm{Mod}$ and $S$-$\mathrm{Mod}$, respectively. Thus there are two exact sequences
$$\xymatrix{R\ar[r]^{\phi}&A_{0}\ar[r]&A_{1}\ar[r]&0}$$
and
$$\xymatrix{S\ar[r]^{\varphi}&B_{0}\ar[r]&B_{1}\ar[r]&0}$$
with $A_{i}\in\mathrm{Add}A$ and $B_{i}\in\mathrm{Add}B$ for $i=0, ~1$, where $\phi$ and $\varphi$ are left $\mathcal{D}_{\sigma_{A}}$-approximation and $\mathcal{D}_{\sigma_{B}}$-approximation, respectively.
Since $F$ is right exact, then there are the following two exact sequences in the comma category.
$$\xymatrix{\binom{R}{FR}\ar[r]&\binom{A_{0}}{FA_{0}}\ar[r]&\binom{A_{1}}{FA_{1}}\ar[r]&0}$$
and
$$\xymatrix{\binom{0}{S}\ar[r]&\binom{0}{B_{0}}\ar[r]&\binom{0}{B_{1}}\ar[r]&0}$$
We consider the following exact sequence in comma category
$$\xymatrix{\binom{R}{FR}\bigoplus\binom{0}{S}\ar[r]^{a}&\binom{A_{0}}{FA_{0}}\bigoplus\binom{0}{B_{0}}\ar[r]&\binom{A_{1}}{FA_{1}}\bigoplus\binom{0}{B_{1}}\ar[r]&0}$$
where $a=\begin{pmatrix}\binom{\phi}{F\phi}&0\\0&\binom{0}{\varphi}\end{pmatrix}$. We claim that the morphism $a$ is left $\mathcal{D}_{\sigma}$-approximation.
For any $\binom{X}{Y}_{h}\in\mathcal{D}_{\sigma}$ and the morphism $(\binom{f}{g},~\binom{0}{x})\in\mathrm{Hom}(\binom{R}{FR}\bigoplus\binom{0}{S},~\binom{X}{Y})$, we need to
find a morphism $\xi$ such that the following diagram is commutative.
$$\xymatrix{
\binom{X}{Y}_{h}&\\
\binom{R}{FR}\bigoplus\binom{0}{S}\ar[u]^{(\binom{f}{g},~\binom{0}{x})}\ar[r]^{a}&\binom{A_{0}}{FA_{0}}\bigoplus\binom{0}{B_{0}}\ar[r]\ar@{.>}[ul]_{\xi}&\binom{A_{1}}{FA_{1}}\bigoplus\binom{0}{B_{1}}\ar[r]&0
}$$
By Lemma \ref{lem1}, $X\in \mathcal{D}_{\sigma_{A}}$ and $Y\in \mathcal{D}_{\sigma_{B}}$, and then there are two morphisms $f_{1}$ and $x_{1}$ such that the following two diagrams are commutative.
$$\xymatrix{
X&\\
R\ar[r]^{\phi}\ar[u]^{f}&A_{0}\ar@{.>}[ul]_{\exists ~f_{1}}\ar[r]&A_{1}\ar[r]&0
}~~~~~~
\xymatrix{
Y&\\
S\ar[r]^{\varphi}\ar[u]^{x}&B_{0}\ar@{.>}[ul]_{\exists ~x_{1}}\ar[r]&B_{1}\ar[r]&0
}$$
We conclude that this morphism $\xi=(\binom{f_{1}}{y},~\binom{0}{x_{1}})$ is what we want, where $y=h\cdot Ff_{1}$.
In fact, we only need to prove that the following two diagrams are commutative, where the first diagram is obvious.
$$\xymatrix{
&0\ar[dr]\ar@{.>}[dd]&\\
0\ar[ur]\ar[dd]\ar[rr]&&FX\ar[dd]\\
&B_{1}\ar@{.>}[dr]^{x_{1}}&\\
S\ar@{.>}[ur]^{\psi}\ar[rr]_{x}&&Y
}~~~~~~~~~
\xymatrix{
&FA_{1}\ar[dr]^{Ff_{1}}\ar@{.>}[dd]&\\
FR\ar[ur]^{F\phi}\ar[dd]\ar[rr]_{Ff}&&FX\ar[dd]^{h}\\
&FA_{1}\ar@{.>}[dr]^{y}&\\
FR\ar@{.>}[ur]^{F\phi}\ar[rr]_{g}&&Y
}$$
Note that $g=h\cdot Ff$ by the definition of $\binom{f}{g}$. Since $y\cdot F\phi=(h\cdot Ff_{1})\cdot F\phi=h\cdot (Ff_{1}\cdot F\phi)=h\cdot F(f_{1}\cdot \phi)=h\cdot Ff=g$,
the second diagram is also commutative. Note that there exists a surjective morphism $b$: $\binom{R^{(X)}}{FR^{(X)}}\bigoplus\binom{0}{S^{(X)}}\longrightarrow \binom{X}{Y}_{h}$ for some set $X$.
According to the above discussion, we can find a morphism $c$ satisfying $b=ca^{(X)}$, and then $c$ is surjective. i.e., $\binom{X}{Y}_{h}\in\mathrm{Gen}(\binom{A}{FA}\bigoplus\binom{0}{B})$.
Then $\binom{A}{FA}\bigoplus\binom{0}{B}$ is silting with respect to $\sigma$.
\ \hfill $\Box$

\vskip 10pt

Let the functor $G$: $\mathrm{Mod}$-$S\longrightarrow\mathrm{Mod}$-$R$ be right exact and covariant.
Similarly, we have the following result in the right comma category $(G$,~$\mathrm{Mod}$-$R$).

\bg{Th}\rm{}\label{3}
Let $X\in \mathrm{Mod}$-$R$, $Y\in \mathrm{Mod}$-$S$ and the functor $G$ commute with direct sums.
Then $(GY\bigoplus X,~Y)_{\binom{1}{0}}$ is silting if and only if
$X$ is silting,
$Y$ is silting and
$GY\in \mathrm{Gen}X$.
\ed{Th}

By Proposition \ref{prop1} and Theorem \ref{th}, we can get the following corollary immediately.

\bg{Cor}\rm{}
Let $A\in R$-$\mathrm{Mod}$, $B\in S$-$\mathrm{Mod}$ and the functor $F$ commute with direct sums.
We have that the following results:

$(1)$ $\binom{A}{FA}\bigoplus\binom{0}{S}$ is silting (resp., partial silting) if and only if $A$ is silting (resp., partial silting).

$(2)$ $\binom{0}{B}$ is is silting (resp., partial silting) if and only if $B$ is silting (resp., partial silting).
\ed{Cor}

\bg{Exm}\rm{}
Let $T$ be given by the quiver
$$
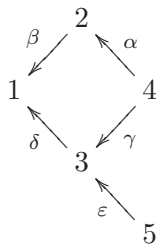
\xymatrix@=0.5cm{
&2\ar[dl]_{\beta}\\
1& &4 \ar[dl]^{\gamma}\ar[ul]_{\alpha} \\
& 3\ar[ul]^{\delta}\\
&&5\ar[ul]^{\varepsilon}
}$$
bounded by $\alpha\beta=\gamma\delta$, $\varepsilon\delta=0$. The Auslander-Reiten quiver of $T$ is as follows.
\ed{Exm}

$$\xymatrix{    &&&&{\smallmatrix5\\3\endsmallmatrix}\ar[dr]& &{\smallmatrix4\\2\endsmallmatrix}\ar[dr]&\\
& {\smallmatrix2\\1\endsmallmatrix}\ar[dr] &   &{\smallmatrix3\endsmallmatrix} \ar[dr]\ar[ur] &
& {\smallmatrix4&5\\2&3\endsmallmatrix}\ar[dr]\ar[ur]   &  &4 \\
1 \ar[dr]\ar[ur]&    & {\smallmatrix2&&3\\&1&\endsmallmatrix}\ar[dr]\ar[ur]\ar[r]  &
{\smallmatrix&4&\\2&&3\\&1&\endsmallmatrix}\ar[r]  &  {\smallmatrix&4&\\2&&3\endsmallmatrix}\ar[ur]\ar[dr]    &
& {\smallmatrix4&&5\\&3&\endsmallmatrix}\ar[ur]\ar[dr] & \\
&  {\smallmatrix3\\1\endsmallmatrix}\ar[ur]  && {\smallmatrix2\endsmallmatrix}\ar[ur]      &
&{\smallmatrix4\\3\endsmallmatrix}\ar[ur] &    &{\smallmatrix5\endsmallmatrix}
}$$

Let $e=e_{1}+e_{3}+e_{5}$, then
$$T\cong\left(\begin{array}{cc}eTe & 0 \\(1-e)Te &(1-e)T(1-e) \\\end{array}\right).$$
Take $R:=eTe(\cong$ $k(\xymatrix{5\ar[r]^{\varepsilon}&3\ar[r]^{\delta}&1})$ with $\varepsilon\delta=0$),
$S:=(1-e)T(1-e)$ ($\cong$$k(\xymatrix{4\ar[r]^{\alpha}&2})$) and $_{S}U_{R}:=(1-e)Te$.
Then the right comma category $(-\bigotimes_{S} U$,~$\mathrm{Mod}$-$R$) is isomorphic to $\mathrm{Mod}$-$T$ \cite{MN}.
Note that $M_{R}={\smallmatrix5\\3\\1\endsmallmatrix}\oplus{\smallmatrix3\\1\endsmallmatrix}\oplus{\smallmatrix3\endsmallmatrix}$
and $N_{S}={\smallmatrix4\\2\endsmallmatrix}\oplus{\smallmatrix4\endsmallmatrix}$ are APR-tilting modules \cite{ASS}.
It is easy to verify that ${\smallmatrix4\\2\endsmallmatrix}\otimes_{S} U={\smallmatrix3\\1\endsmallmatrix}$,
${\smallmatrix2\endsmallmatrix}\otimes_{S} U={\smallmatrix1\endsmallmatrix}$, and ${\smallmatrix4\endsmallmatrix}\otimes_{S} U={\smallmatrix3\endsmallmatrix}$.
Then $N\otimes_{S} U={\smallmatrix3\endsmallmatrix}\oplus{\smallmatrix3\\1\endsmallmatrix}\in\mathrm{Gen}M$.
By the proposition 3.13 in \cite{AMV}, both $M_{R}$ and $N_{S}$ are silting moduels,
and then $(N\otimes_{S} U,~N)\bigoplus(M,~0)$ is silting in $\mathrm{Mod}$-$T$ by Theorem \ref{3}.

\section{Torsion pairs}

Let $U$ is an $(S,~R)$-bimodule. Take $F=U\bigotimes_{R}-$, then the comma category $(F,~S$-$\mathrm{Mod})$ is isomorphic to the category $T$-Mod \cite{MN}, where
$T=\left(\begin{array}{cc}R & 0 \\U & S \\\end{array}\right)$. In this subsection, we mainly consider the special comma category $(U\bigotimes_{R}-,~S$-$\mathrm{Mod})\cong T$-$\mathrm{Mod}$.
In this case, for any object $\binom{A}{B}$ is in $(U\bigotimes_{R}-,~S$-$\mathrm{Mod})$, we have that $A\in R$-Mod and $B\in S$-Mod.

Firstly, we study the following three functors on the product category and the comma category and three classes, which are very useful.

\bg{Def}\rm{}\label{}%
$(1)$ $\textbf{p}$: $R$-$\mathrm{Mod}\times S$-$\mathrm{Mod}\longrightarrow$ $T$-$\mathrm{Mod}$ is defined as follows: for each objects $(A,~B)$ of $R$-$\mathrm{Mod}\times S$-$\mathrm{Mod}$,
let $\textbf{p}(A,~B)=\binom{A}{(U\bigotimes_{R}A)\bigoplus B}_{\binom{1}{0}}$ with the obvious map and for any morphism $(f,~g)$ in $R$-$\mathrm{Mod}\times S$-$\mathrm{Mod}$,
let $\textbf{p}(f,~g)=\binom{f}{(U\bigotimes_{R}f)\bigoplus g}$. It is easy to see that $\textbf{p}(A,~B)=\textbf{p}(A,~0)\bigoplus\textbf{p}(0,~B)$.

$(2)$ $\textbf{q}$: $T$-$\mathrm{Mod}$ $\longrightarrow R$-$\mathrm{Mod}\times S$-$\mathrm{Mod}$ is defined as follows: for each objects $\binom{A}{B}$ of $T$-$\mathrm{Mod}$,
let $\textbf{q}\binom{A}{B}=(A,~B)$ with the obvious map and for any morphism $\binom{f}{g}$ in $T$-$\mathrm{Mod}$,
let $\textbf{q}\binom{f}{g}=(f,~g)$.

$(3)$ $\textbf{h}$: $R$-$\mathrm{Mod}\times S$-$\mathrm{Mod}\longrightarrow$ $T$-$\mathrm{Mod}$ is defined as follows:
for each objects $(A,~B)$ of $R$-$\mathrm{Mod}\times S$-$\mathrm{Mod}$,
let $\textbf{h}(A,~B)=\binom{A\bigoplus \mathrm{Hom}_{S}(U,B)}{ B}$ with the obvious map and for any morphism $(f,~g)$ in $R$-$\mathrm{Mod}\times S$-$\mathrm{Mod}$,
let $\textbf{h}(f,~g)=\binom{f\oplus \mathrm{Hom}_{S}(U,g)}{g}$.

It is easy to verify that $\textbf{p}$ is a left adjoint of $\textbf{q}$ and $\textbf{q}$ is a left adjoint of $\textbf{h}$.
\ed{Def}

The functor $\textbf{p}$ was introduced by Mitchell \cite{M}, which was a particular case of the additive
left Kan extension functor. Recently, the functor $\textbf{p}$ was used by Enochs \cite{ECT} and Hu \cite{HZ}. The functors $\textbf{q}$ and $\textbf{h}$ was used by Mao \cite{MLX}.

Let $\mathcal{C}$ be a class of $R$-$\mathrm{Mod}$ and $\mathcal{D}$ be a class of $S$-$\mathrm{Mod}$. We will denote by
$\mathfrak{U}^{\mathcal{C}}_{\mathcal{D}}$ the class of objects $\{\binom{C}{D}$: $C\in\mathcal{C}$ and $D\in\mathcal{D}\}$. Denoted by
$\mathfrak{B}^{\mathcal{C}}_{\mathcal{D}}$ the class of objects $\{\binom{C}{D}_{\varphi}$: $C\in\mathcal{C}$ and $D/\mathrm{Im}(\varphi)\in\mathcal{D}$, $\varphi$ is a monomorphism$\}$.
Denoted by $\mathfrak{J}^{\mathcal{C}}_{\mathcal{D}}$ the class of objects $\{\binom{C}{D}_{\varphi}$: $\mathrm{ker}\widetilde{\varphi}\in\mathcal{C}$ and
$D\in\mathcal{D}$, $\widetilde{\varphi}$ is an epimorphism$\}$, where $\widetilde{\varphi}$ is the morphism from $A$ to $\mathrm{Hom}_{S}(U,~B)$ given by
$\widetilde{\varphi}(x)(u)=\varphi(u\otimes x)$ for each $u\in U$ and $x\in A$.

Next, we will give some isomorphisms on $\mathrm{Hom}$ set in $T$-$\mathrm{Mod}$, which are very important in the next work.

\bg{Lem}\rm{}\label{lem4}%
Let $\binom{A}{B}$ and $\binom{C}{D}_{\varphi}$ be in $T$-$\mathrm{Mod}$. Then the following statements hold.

$(1)$ $\mathrm{Hom}_{T}(\binom{A}{B},~\binom{C}{0})\cong \mathrm{Hom}_{R}(A,~C)$.

$(2)$ $\mathrm{Hom}_{T}(\binom{0}{B},~\binom{C}{D})\cong \mathrm{Hom}_{S}(B,~D)$.

$(3)$ $\mathrm{Hom}_{T}(\binom{A}{U\bigotimes_{R}A},~\binom{C}{D})\cong \mathrm{Hom}_{R}(A,~C)$.

$(4)$ $\mathrm{Hom}_{T}(\binom{A}{B},~\binom{\mathrm{Hom}_{S}(U,~D)}{D})\cong \mathrm{Hom}_{S}(B,~D)$.

$(5)$ $\mathrm{Hom}_{T}(\binom{R}{0},~\binom{C}{D})\cong \mathrm{ker}\widetilde{\varphi}$.
\ed{Lem}

\Pf. (1) Since the pair (\textbf{q},~\textbf{h}) is an adjoint pair, we have that $\mathrm{Hom}_{T}(\binom{A}{B},~\binom{C}{0})\cong \mathrm{Hom}_{T}(\binom{A}{B},~\textbf{h}(C,~0))$
$\cong\mathrm{Hom}_{R-\mathrm{Mod}\times S-\mathrm{Mod}}(\textbf{q}\binom{A}{B},~(C,~0))$ $\cong\mathrm{Hom}_{R-\mathrm{Mod}\times S-\mathrm{Mod}}((A,~B),~(C,~0))$ $\cong\mathrm{Hom}_{R}(A,~C)$.

(2)--(4) Since $(\textbf{p},~\textbf{q})$, $(\textbf{q},~\textbf{h})$ are adjoint pairs, these proofs are similar to (1).

(5) holds by Lemma 4.1(2) in \cite{MLX}.
\ \hfill $\Box$

\vskip 10pt

Analogously, we have the category Mod-$T$ whose objects are triples $(X,~Y)_{\phi}$, where $X\in$ Mod-$R$, $Y\in$ Mod-$S$ and $\phi$: $Y\bigotimes_{B} U\longrightarrow X$ is an $R$-module morphism,
and whose morphisms from $(X_{1},~Y_{1})_{\phi_{1}}$ to $(X_{2},~Y_{2})_{\phi_{2}}$ are pairs $(g_{1},~g_{2})$ such that $g_{1}\in\mathrm{Hom}_{R}(X_{1},~X_{2})$, $g_{2}\in\mathrm{Hom}_{S}(Y_{1},~Y_{2})$
and $\phi_{2}(g_{2}\otimes 1)=g_{1}\phi_{1}$.

By \cite[Proposition 3.6.1]{KT}, there is an isomorphism of abelian groups
$$(X,~Y)_{\phi}\bigotimes_{T}\binom{A}{B}_{\varphi}\cong((X\bigotimes_{R} A)\bigoplus(Y\bigotimes_{S} B))/H$$
where the subgroup $H$ is generated by all elements of the form $(\phi(y\otimes u))\otimes a-y\otimes \varphi(u\otimes a)$
with $y\in Y$, $u\in U$ and $a\in A$. It follows from the definition of tensor products that the following conclusions hold.

\bg{Pro}\rm{}\label{prop2}%
Let $(A,~B)$ be a right $T$-module and $\binom{C}{D}_{\varphi}$ be a left $T$-module. Then the following statements hold.

$(1)$ $(A,~0)\bigotimes_{T}\binom{C}{D}\cong A\bigotimes_{R}C$.

$(2)$ $(A,~B)\bigotimes_{T}\binom{0}{D}\cong B\bigotimes_{S}D$.

$(3)$ $(A,~B)\bigotimes_{T}\binom{C}{U\bigotimes_{R}C}\cong A\bigotimes_{R}C$.

$(4)$ $(B\bigotimes_{S}U,~B)\bigotimes_{T}\binom{C}{D}\cong B\bigotimes_{S}D$.

$(5)$ \cite[Lemma 4.1(1)]{MLX} $(0,~B)\bigotimes_{T}\binom{C}{D}\cong D/\mathrm{Im}\varphi$.
\ed{Pro}

The following proposition gives a very important property of the two classes $\mathfrak{U}^{\mathcal{C}}_{\mathcal{D}}$ and $\mathfrak{B}^{\mathcal{C}}_{\mathcal{D}}$.

\bg{Pro}\rm{}\label{prop3}
Let $\mathcal{C}$ be a class of left $R$-modules and $\mathcal{D}$ be a class of left $S$-modules. Then the following statements hold.

$(1)$ $(\mathfrak{B}^{\mathcal{C}}_{\mathcal{D}})^{\perp_{0}}=\mathfrak{U}^{\mathcal{C}^{\perp_{0}}}_{\mathcal{D}^{\perp_{0}}}$.

$(2)$ $\mathfrak{B}^{^{\perp_{0}}\mathcal{C}}_{^{\perp_{0}}\mathcal{D}}\subseteq$ $^{\perp_{0}}(\mathfrak{U}^{\mathcal{C}}_{\mathcal{D}})$.
The converses hold if $S^{+}\in\mathcal{D}$, $_{R}U$ is flat and $U_{S}$ is projective.
%
%
\ed{Pro}

\Pf. (1) For any $\binom{A}{B}\in(\mathfrak{B}^{\mathcal{C}}_{\mathcal{D}})^{\perp_{0}}$, $C\in\mathcal{C}$ and $D\in \mathcal{D}$.
Since $\binom{C}{U\bigotimes_{R}C}\in\mathfrak{B}^{\mathcal{C}}_{\mathcal{D}}$, we have that
$\binom{C}{U\bigotimes_{R}C}\in$ $\mathrm{Hom}_{R}(C,~A)\cong \mathrm{Hom}(\binom{C}{U\bigotimes_{R}C},~\binom{A}{B})=0$ by Lemma \ref{lem4}(3).
Since $\binom{0}{D}\in\mathfrak{B}^{\mathcal{C}}_{\mathcal{D}}$, we have that $\mathrm{Hom}_{S}(D,~B)\cong \mathrm{Hom}(\binom{0}{D},~\binom{A}{B})=0$ by Lemma \ref{lem4}(2).
i.e., $(\mathfrak{B}^{\mathcal{C}}_{\mathcal{D}})^{\perp_{0}}\subseteq\mathfrak{U}^{\mathcal{C}^{\perp_{0}}}_{\mathcal{D}^{\perp_{0}}}$.

On the contrary, For any $\binom{M}{N}\in\mathfrak{U}^{\mathcal{C}^{\perp_{0}}}_{\mathcal{D}^{\perp_{0}}}$ and $\binom{X}{Y}\in\mathfrak{B}^{\mathcal{C}}_{\mathcal{D}}$.
We can obtain the following exact sequence in $T$-$\mathrm{Mod}$ by the definition of $\mathfrak{B}^{\mathcal{C}}_{\mathcal{D}}$.
$$\xymatrix{0\ar[r]&\binom{X}{U\bigotimes_{R}X}\ar[r]&\binom{X}{Y}\ar[r]&\binom{0}{Y/U\bigotimes_{R}X}\ar[r]&0}.$$
Applying the functor $\mathrm{Hom}_{T}(-,~\binom{M}{N})$ to the above exact sequence, we can get the following exact sequence
$$\xymatrix{0\ar[r]&\mathrm{Hom}_{T}(\binom{0}{Y/U\bigotimes_{R}X},~\binom{M}{N})\ar[r]&\mathrm{Hom}_{T}(\binom{X}{Y},~\binom{M}{N})\ar[r]
&\mathrm{Hom}_{T}(\binom{X}{U\bigotimes_{R}X},~\binom{M}{N})}.$$
Note that $X\in C$ and $Y/U\bigotimes_{R}X\in D$ by the definition of $\mathfrak{B}^{\mathcal{C}}_{\mathcal{D}}$.
We have that $\mathrm{Hom}_{T}(\binom{0}{Y/U\bigotimes_{R}X},~\binom{M}{N})$ $\cong\mathrm{Hom}_{S}(Y/U\bigotimes_{R}X,~N)=0$ and
$\mathrm{Hom}_{T}(\binom{X}{U\bigotimes_{R}X},~\binom{M}{N})\cong\mathrm{Hom}_{R}(X,~M)$ by Lemma \ref{lem4} (2) and (3).
So we get that $\mathrm{Hom}_{T}(\binom{X}{Y},~\binom{M}{N})=0$ from the above exact sequence.
i.e., $\mathfrak{U}^{\mathcal{C}^{\perp_{0}}}_{\mathcal{D}^{\perp_{0}}}\subseteq(\mathfrak{B}^{\mathcal{C}}_{\mathcal{D}})^{\perp_{0}}$.

(2) For any $\binom{G}{H}\in\mathfrak{B}^{^{\perp_{0}}\mathcal{C}}_{^{\perp_{0}}\mathcal{D}}$, $P\in\mathcal{C}$ and $Q\in \mathcal{D}$.
We can obtain the following exact sequence in $T$-$\mathrm{Mod}$ by the definition of $\mathfrak{B}^{^{\perp_{0}}\mathcal{C}}_{^{\perp_{0}}\mathcal{D}}$.
$$\xymatrix{0\ar[r]&\binom{G}{U\bigotimes_{R}G}\ar[r]&\binom{G}{H}\ar[r]&\binom{0}{H/U\bigotimes_{R}G}\ar[r]&0}$$
with $G\in$ $^{\perp_{0}}\mathcal{C}$ and $H/U\bigotimes_{R}G\in$ $^{\perp_{0}}\mathcal{D}$.
Applying the functor $\mathrm{Hom}_{T}(-,~\binom{P}{Q})$ to the above exact sequence, we can obtain the following exact sequence.
$$\xymatrix{0\ar[r]&\mathrm{Hom}_{T}(\binom{0}{H/U\bigotimes_{R}G},~\binom{P}{Q})\ar[r]&\mathrm{Hom}_{T}(\binom{G}{H},~\binom{P}{Q})\ar[r]&\mathrm{Hom}_{T}(\binom{G}{U\bigotimes_{R}G},~\binom{P}{Q})}.$$
Note that $\mathrm{Hom}_{T}(\binom{0}{H/U\bigotimes_{R}G},~\binom{P}{Q})\cong\mathrm{Hom}_{S}(H/U\bigotimes_{R}G,~Q)=0$ and
$\mathrm{Hom}_{T}(\binom{G}{U\bigotimes_{R}G},~\binom{P}{Q})\cong\mathrm{Hom}_{R}(G,~P)=0$ by the lemma \ref{lem4} (2) and (3).
Consequently, we get that $\mathrm{Hom}_{T}(\binom{G}{H},~\binom{P}{Q})=0$ from the above exact sequence.
i.e., $\mathfrak{B}^{^{\perp_{0}}\mathcal{C}}_{^{\perp_{0}}\mathcal{D}}\subseteq$ $^{\perp_{0}}(\mathfrak{U}^{\mathcal{C}}_{\mathcal{D}})$.

Conversely, $\binom{I}{J}\in$ $^{\perp_{0}}(\mathfrak{U}^{\mathcal{C}}_{\mathcal{D}})$, $C\in\mathcal{C}$ and $D\in \mathcal{D}$.
By Lemma \ref{lem4} (1), we have that $\mathrm{Hom}_{R}(I,~C)\cong\mathrm{Hom}_{T}(\binom{I}{J},~\binom{C}{0})=0$, and then $I\in$ $^{\perp_{0}}\mathcal{C}$.

Note that $((0~,S)\bigotimes_{T}-,~\mathrm{Hom}_{S}((0~,S),-)))$ is an adjoint pair. From the following isomorphisms, we have that $J/(U\bigotimes_{R}I)\in$ $^{\perp_{0}}\mathcal{D}$.
\begin{align*}
\mathrm{Hom}_{S}(J/(U\bigotimes_{R}I),~D)& \cong  \mathrm{Hom}_{S}((0~,S)\bigotimes_{T}\binom{I}{J},~D), \rm{}by ~the~ definition ~of ~tensor ~pruducts.\\
&\cong\mathrm{Hom}_{T}(\binom{I}{J},~\mathrm{Hom}_{S}((0~,S),~D)), \rm{}by ~the~ definition ~of ~adjoint ~pairs.\\
& \cong\mathrm{Hom}_{T}(\binom{I}{J},~\binom{0}{D})\\
& =0
\end{align*}

We consider the following exact sequence in category Mod-$T$.
$$\xymatrix{0\ar[r]&(U,~0)\ar[r]&(U,~S)\ar[r]&(0,~S)\ar[r]&0}.$$
Applying the functor $-\bigotimes_{T}\binom{I}{J}$ to the above exact sequence, we can obtain the following exact sequence.
$$\xymatrix{\mathrm{Tor}_{1}((0,~S),~\binom{I}{J})\ar[r]&(U,~0)\bigotimes_{T}\binom{I}{J}\ar[r]&(U,~S)\bigotimes_{T}\binom{I}{J}\ar[r]&(0,~S)\bigotimes_{T}\binom{I}{J}\ar[r]&0}.$$
By the proposition \ref{prop2} (1) and (4), we can get that $(U,~0)\bigotimes_{T}\binom{I}{J}\cong U\bigotimes_{R} I$ and
$(U,~S)\bigotimes_{T}\binom{I}{J}\cong(S\bigotimes_{S}U,~S)\bigotimes_{T}\binom{I}{J}\cong S\bigotimes_{S}J\cong J$.
Note $\mathrm{Tor}_{1}((0,~S),~\binom{I}{J})^{+}\cong \mathrm{Ext}_{1}(\binom{I}{J},~\binom{0}{S^{+}})$.
Since $S^{+}$ is injective, the injective dimension of $\binom{0}{S^{+}}$ is at most 1 by the corollary 3.3 (2) in \cite{MLX}.
It is easy to verify from $\mathrm{Hom}_{T}(\binom{I}{J},~\binom{0}{S^{+}})=0$ that $\mathrm{Ext}_{1}(\binom{I}{J},~\binom{0}{S^{+}})=0$. i.e., $\mathrm{Tor}_{1}((0,~S)\bigotimes_{T}\binom{I}{J})=0$.
So the morphism $U\bigotimes_{R} I\longrightarrow J$ is an monomorphism. i.e.,
$^{\perp_{0}}(\mathfrak{U}^{\mathcal{C}}_{\mathcal{D}})\subseteq$ $\mathfrak{B}^{^{\perp_{0}}\mathcal{C}}_{^{\perp_{0}}\mathcal{D}}$.
\ \hfill $\Box$

\vskip 10pt
From the proof of (2) in the above proposition, we can use condition "$\mathrm{Tor}_{1}((0,~S)$, $^{\perp_{0}}(\mathfrak{U}^{\mathcal{C}}_{\mathcal{D}}))=0$
instead of condition "$S^{+}\in\mathcal{D}$, $_{R}U$ is flat and $U_{S}$ is projective". Unfortunately, we don't know whether the conclusion hold after deleting the condition.
Next, we will give the theorem \ref{2} (1) and its proof.

\bg{Th}\rm{}\label{th1}
Let $\mathcal{C}_{1}$ and $\mathcal{C}_{2}$ be two subcategories of left $R$-modules and $\mathcal{D}_{1}$ and $\mathcal{D}_{2}$ be two subcategories of left $S$-modules.
If $S^{+}\in\mathcal{D}$, $_{R}U$ is flat and $U_{S}$ is projective, then $(\mathcal{C}_{1},~\mathcal{C}_{2})$ and $(\mathcal{D}_{1},~\mathcal{D}_{2})$ are torsion pairs if and only if
$(\mathfrak{B}^{\mathcal{C}_{1}}_{\mathcal{D}_{1}},~\mathfrak{U}^{\mathcal{C}_{2}}_{\mathcal{D}_{2}})$ is a torsion pair.
\ed{Th}

\Pf. ($\Rightarrow$) Since $(\mathcal{C}_{1},~\mathcal{C}_{2})$ and $(\mathcal{D}_{1},~\mathcal{D}_{2})$ are torsion pairs, we have that
$(\mathfrak{B}^{\mathcal{C}_{1}}_{\mathcal{D}_{1}})^{\perp_{0}}=\mathfrak{U}^{\mathcal{C}_{1}^{\perp_{0}}}_{\mathcal{D}_{1}^{\perp_{0}}}$
$=\mathfrak{U}^{\mathcal{C}_{2}}_{\mathcal{D}_{2}}$ by the proposition \ref{prop3} (1).
And $^{\perp_{0}}(\mathfrak{U}^{\mathcal{C}_{2}}_{\mathcal{D}_{2}})=\mathfrak{B}^{^{\perp_{0}}\mathcal{C}_{2}}_{^{\perp_{0}}\mathcal{D}_{2}}=\mathfrak{B}^{\mathcal{C}_{1}}_{\mathcal{D}_{1}}$
by Proposition \ref{prop3} (2).

($\Leftarrow$) (1) For any $A\in\mathcal{C}_{1}$ and $B\in\mathcal{C}_{2}$, then $\binom{A}{U\bigotimes_{R}A}\in\mathfrak{B}^{\mathcal{C}_{1}}_{\mathcal{D}_{1}}$
and $\binom{B}{0}\in\mathfrak{U}^{\mathcal{C}_{2}}_{\mathcal{D}_{2}}$.
Since $(\mathfrak{B}^{\mathcal{C}_{1}}_{\mathcal{D}_{1}},~\mathfrak{U}^{\mathcal{C}_{2}}_{\mathcal{D}_{2}})$ is a torsion pair,
we have that $\mathrm{Hom}_{R}(A,~B)\cong \mathrm{Hom}_{T}(\binom{A}{U\bigotimes_{R}A},~\binom{B}{0})=0$ by Lemma \ref{lem4} (3).
i.e., $\mathrm{Hom}_{R}(\mathcal{C}_{1},~\mathcal{C}_{2})=0$

For any $X\in\mathcal{C}_{1}^{\perp_{0}}$, then $\binom{X}{0}\in\mathfrak{U}^{\mathcal{C}_{1}^{\perp_{0}}}_{\mathcal{D}_{1}^{\perp_{0}}}=
(\mathfrak{B}^{\mathcal{C}_{1}}_{\mathcal{D}_{1}})^{\perp_{0}}=\mathfrak{U}^{\mathcal{C}_{2}}_{\mathcal{D}_{2}}$ by Proposition \ref{prop3} (1). i.e., $X\in\mathcal{C}_{2}$.

For any $Y\in$ $^{\perp_{0}}\mathcal{C}_{2}$, then $\binom{Y}{U\bigotimes_{R}Y}\in\mathfrak{B}^{^{\perp_{0}}\mathcal{C}_{2}}_{^{\perp_{0}}\mathcal{D}_{2}}=$
$^{\perp_{0}}(\mathfrak{U}^{\mathcal{C}_{2}}_{\mathcal{D}_{2}})=\mathfrak{B}^{\mathcal{C}_{1}}_{\mathcal{D}_{1}}$ by Proposition \ref{prop3} (2). i.e., $Y\in\mathcal{C}_{1}$.
So $(\mathcal{C}_{1},~\mathcal{C}_{2})$ is a torsion pair.

(2) For any $M\in\mathcal{D}_{1}$ and $N\in\mathcal{D}_{2}$, then $\binom{0}{M}\in\mathfrak{B}^{\mathcal{C}_{1}}_{\mathcal{D}_{1}}$
and $\binom{0}{N}\in\mathfrak{U}^{\mathcal{C}_{2}}_{\mathcal{D}_{2}}$. By Lemma \ref{lem4} (2), we have that
$\mathrm{Hom}_{S}(M,~N)\cong \mathrm{Hom}_{T}(\binom{0}{M},~\binom{0}{N})=0$
since $(\mathfrak{B}^{\mathcal{C}_{1}}_{\mathcal{D}_{1}},~\mathfrak{U}^{\mathcal{C}_{2}}_{\mathcal{D}_{2}})$ is a torsion pair.

For any $P\in\mathcal{D}_{1}^{\perp_{0}}$, then $\binom{0}{P}\in\mathfrak{U}^{\mathcal{C}_{1}^{\perp_{0}}}_{\mathcal{D}_{1}^{\perp_{0}}}=
(\mathfrak{B}^{\mathcal{C}_{1}}_{\mathcal{D}_{1}})^{\perp_{0}}=\mathfrak{U}^{\mathcal{C}_{2}}_{\mathcal{D}_{2}}$ by Proposition \ref{prop3} (1). i.e., $P\in\mathcal{D}_{2}$.

For any $Q\in$ $^{\perp_{0}}\mathcal{D}_{2}$, then $\binom{0}{Q}\in\mathfrak{B}^{^{\perp_{0}}\mathcal{C}_{2}}_{^{\perp_{0}}\mathcal{D}_{2}}=$
$^{\perp_{0}}(\mathfrak{U}^{\mathcal{C}_{2}}_{\mathcal{D}_{2}})=\mathfrak{B}^{\mathcal{C}_{1}}_{\mathcal{D}_{1}}$ by Proposition \ref{prop3} (2). i.e., $Q\in\mathcal{D}_{1}$.
So $(\mathcal{D}_{1},~\mathcal{D}_{2})$ is a torsion pair.
\ \hfill $\Box$

\vskip 10pt
The following proposition gives a connection between $\mathfrak{U}^{\mathcal{C}}_{\mathcal{D}}$ and $\mathfrak{J}^{\mathcal{C}}_{\mathcal{D}}$.

\bg{Pro}\rm{}\label{prop4}
Let $\mathcal{C}$ be a class of left $R$-modules and $\mathcal{D}$ be a class of left $S$-modules. Then the following statements hold.
%
%

$(1)$ $^{\perp_{0}}(\mathfrak{J}^{\mathcal{C}}_{\mathcal{D}})=\mathfrak{U}^{^{\perp_{0}}\mathcal{C}}_{^{\perp_{0}}\mathcal{D}}$.

$(2)$ $\mathfrak{J}^{\mathcal{C}^{\perp_{0}}}_{\mathcal{D}^{\perp_{0}}}\subseteq(\mathfrak{U}^{\mathcal{C}}_{\mathcal{D}})^{\perp_{0}}$.
The converses hold if $R\in\mathcal{C}$, $_{R}U$ is flat and $U_{S}$ is projective.
\ed{Pro}

\Pf. (1) ($\Rightarrow$) For any $\binom{X}{Y}\in$ $^{\perp_{0}}(\mathfrak{J}^{\mathcal{C}}_{\mathcal{D}})$, $C\in\mathcal{C}$ and $D\in \mathcal{D}$.
Note that $\binom{C}{0}\in\mathfrak{J}^{\mathcal{C}}_{\mathcal{D}}$. By Lemma \ref{lem4} (1), we have that
$\mathrm{Hom}_{R}(X,~C)\cong\mathrm{Hom}_{T}(\binom{X}{Y},~\binom{C}{0})=0$. i.e., $X\in$ $^{\perp_{0}}\mathcal{C}$.

Note that $\binom{\mathrm{Hom}_{S}(U,~D)}{D}\in\mathfrak{J}^{\mathcal{C}}_{\mathcal{D}}$. By Lemma \ref{lem4} (4), we have that
$\mathrm{Hom}_{S}(Y,~D)\cong\mathrm{Hom}_{T}(\binom{X}{Y}$, $\binom{\mathrm{Hom}_{S}(U,~D)}{D})=0$. i.e., $Y\in$ $^{\perp_{0}}\mathcal{D}$.
So $^{\perp_{0}}(\mathfrak{J}^{\mathcal{C}}_{\mathcal{D}})\subseteq\mathfrak{U}^{^{\perp_{0}}\mathcal{C}}_{^{\perp_{0}}\mathcal{D}}$.

($\Leftarrow$) For any $\binom{M}{N}\in\mathfrak{U}^{^{\perp_{0}}\mathcal{C}}_{^{\perp_{0}}\mathcal{D}}$, $\binom{A}{B}_{\varphi}\in\mathfrak{J}^{\mathcal{C}}_{\mathcal{D}}$.
By the definition of $\mathfrak{J}^{\mathcal{C}}_{\mathcal{D}}$, we can obtain the following exact sequence.
$$\xymatrix{0\ar[r]&\binom{\mathrm{ker}\widetilde{\varphi}}{0}\ar[r]&\binom{A}{B}\ar[r]&\binom{\mathrm{Hom}_{S}(U,~B)}{B}\ar[r]&0}$$
Applying the functor $\mathrm{Hom}_{T}(\binom{M}{N},-)$ to the above exact sequence, we can obtain the following exact sequence.
$$\xymatrix{0\ar[r]&\mathrm{Hom}_{T}(\binom{M}{N},~\binom{\mathrm{ker}\widetilde{\varphi}}{0})\ar[r]&\mathrm{Hom}_{T}(\binom{M}{N},~\binom{A}{B})\ar[r]&\mathrm{Hom}_{T}(\binom{M}{N},~\binom{\mathrm{Hom}_{S}(U,~B)}{B})}$$
By Lemma \ref{lem4} (1) and (4), we have that $\mathrm{Hom}_{T}(\binom{M}{N},~\binom{\mathrm{ker}\widetilde{\varphi}}{0})\cong \mathrm{Hom}_{R}(M,~\mathrm{ker}\widetilde{\varphi})=0$
and $\mathrm{Hom}_{T}(\binom{M}{N},~\binom{\mathrm{Hom}_{S}(U,~B)}{B})\cong\mathrm{Hom}_{S}(N,~B)=0$. From the above sequence, $\mathrm{Hom}_{T}(\binom{M}{N},~\binom{A}{B})=0$.
i.e., $\binom{M}{N}\in$ $^{\perp_{0}}(\mathfrak{J}^{\mathcal{C}}_{\mathcal{D}})$.

(2) For any $\binom{G}{H}_{\phi}\in\mathfrak{J}^{\mathcal{C}^{\perp_{0}}}_{\mathcal{D}^{\perp_{0}}}$, $\binom{I}{J}\in\mathfrak{U}^{\mathcal{C}}_{\mathcal{D}}$.
We consider the following exact sequence.
$$\xymatrix{0\ar[r]&\binom{\mathrm{ker}\widetilde{\phi}}{0}\ar[r]&\binom{G}{H}\ar[r]&\binom{\mathrm{Hom}_{S}(U,~H)}{H}\ar[r]&0}$$
Applying the functor $\mathrm{Hom}_{T}(\binom{I}{J},-)$ to the above exact sequence, we can obtain the following exact sequence.
$$\xymatrix{0\ar[r]&\mathrm{Hom}_{T}(\binom{I}{J},~\binom{\mathrm{ker}\widetilde{\phi}}{0})\ar[r]&\mathrm{Hom}_{T}(\binom{I}{J},~\binom{G}{H})\ar[r]&\mathrm{Hom}_{T}(\binom{I}{J},~\binom{\mathrm{Hom}_{S}(U,~H)}{H})}$$
By Lemma \ref{lem4} (1) and (4), we have that $\mathrm{Hom}_{T}(\binom{I}{J},~\binom{\mathrm{ker}\widetilde{\phi}}{0})\cong \mathrm{Hom}_{R}(I,~\mathrm{ker}\widetilde{\varphi})=0$
and $\mathrm{Hom}_{T}(\binom{I}{J},~\binom{\mathrm{Hom}_{S}(U,~H)}{H})\cong\mathrm{Hom}_{S}(J,~H)=0$. From the above sequence, $\mathrm{Hom}_{T}(\binom{I}{J},~\binom{G}{H})=0$.
So we can obtain that $\mathfrak{J}^{\mathcal{C}^{\perp_{0}}}_{\mathcal{D}^{\perp_{0}}}\subseteq(\mathfrak{U}^{\mathcal{C}}_{\mathcal{D}})^{\perp_{0}}$.

Conversely, For any $\binom{P}{Q}_{\psi}\in(\mathfrak{U}^{\mathcal{C}}_{\mathcal{D}})^{\perp_{0}}$, $C\in \mathcal{C}$ and $D\in \mathcal{D}$.

By Lemma \ref{lem4} (5), we get that
$\mathrm{Hom}_{R}(C,~\mathrm{ker}\widetilde{\psi})\cong\mathrm{Hom}_{R}(C,~\mathrm{Hom}_{T}(\binom{R}{0},~\binom{P}{Q}))
\cong\mathrm{Hom}_{T}(\binom{R}{0}\bigotimes_{R}C,~\binom{P}{Q})\cong\mathrm{Hom}_{T}(\binom{C}{0},~\binom{P}{Q})=0$
since $\binom{C}{0}\in\mathfrak{U}^{\mathcal{C}}_{\mathcal{D}}$. i.e., $\mathrm{ker}\widetilde{\psi}\in\mathcal{C}^{\perp_{0}}$.

Note that $\binom{0}{D}\in\mathfrak{U}^{\mathcal{C}}_{\mathcal{D}}$, and then $\mathrm{Hom}_{S}(D,~Q)\cong\mathrm{Hom}_{T}(\binom{0}{D},~\binom{P}{Q})=0$.
i.e., $Q\in\mathcal{D}^{\perp_{0}}$.

We consider the following exact sequence.
$$\xymatrix{0\ar[r]&\binom{0}{U}\ar[r]&\binom{R}{U}\ar[r]&\binom{R}{0}\ar[r]&0}$$
Applying the functor $\mathrm{Hom}_{T}(-,~\binom{P}{Q})$ to the above exact sequence, we can obtain the following exact sequence.
$$\xymatrix{0\ar[r]&\mathrm{Hom}_{T}(\binom{R}{0},~\binom{P}{Q})\ar[r]&\mathrm{Hom}_{T}(\binom{R}{U},~\binom{P}{Q})\ar[r]&\mathrm{Hom}_{T}(\binom{0}{U},~\binom{P}{Q})\ar[r]&\mathrm{Ext}_{1}(\binom{R}{0},~\binom{P}{Q})}$$
The projective dimension of $\binom{R}{0}$ is at most 1 by Corollary 3.3 (1) in \cite{MLX}.
It is easy to verify from $\mathrm{Hom}_{T}(\binom{R}{0},~\binom{P}{Q})=0$ that $\mathrm{Ext}_{1}(\binom{R}{0},~\binom{P}{Q})=0$.
Note that $\mathrm{Hom}_{T}(\binom{R}{U},~\binom{P}{Q})\cong\mathrm{Hom}_{T}(\binom{R}{U\bigotimes_{R}R},~\binom{P}{Q})
\cong\mathrm{Hom}_{R}(R,~P)\cong P$ and $\mathrm{Hom}_{T}(\binom{0}{U},~\binom{P}{Q})\cong\mathrm{Hom}_{S}(U,~Q)$ by Lemma \ref{lem4} (2) and (3).
i.e., $\binom{P}{Q}_{\psi}\in\mathfrak{J}^{\mathcal{C}^{\perp_{0}}}_{\mathcal{D}^{\perp_{0}}}$.
So we have that $(\mathfrak{U}^{\mathcal{C}}_{\mathcal{D}})^{\perp_{0}}\subseteq\mathfrak{J}^{\mathcal{C}^{\perp_{0}}}_{\mathcal{D}^{\perp_{0}}}$.
\ \hfill $\Box$

\vskip 10pt
Similar to Proposition \ref{prop3} (2), we can use condition "$\mathrm{Ext}^{1}(\binom{R}{0},~(\mathfrak{U}^{\mathcal{C}}_{\mathcal{D}})^{\perp_{0}})=0$"
instead of condition "$R\in\mathcal{C}$, $_{R}U$ is flat and $U_{S}$ is projective". We still don't know whether the conclusion hold or not after deleting the condition.
Finally, we will give Theorem \ref{2} (2) and its proof.

\bg{Th}\rm{}\label{th2}
Let $\mathcal{C}_{1}$ and $\mathcal{C}_{2}$ be two subcategories of left $R$-modules, $\mathcal{D}_{1}$ and $\mathcal{D}_{2}$ be two subcategories of left $S$-modules.
If $R\in\mathcal{C}_{1}$, $_{R}U$ is flat and $U_{S}$ is projective, then $(\mathcal{C}_{1},~\mathcal{C}_{2})$ and $(\mathcal{D}_{1},~\mathcal{D}_{2})$ are torsion pairs if and only if
$(\mathfrak{U}^{\mathcal{C}_{1}}_{\mathcal{D}_{1}},~\mathfrak{J}^{\mathcal{C}_{2}}_{\mathcal{D}_{2}})$ is a torsion pair.
\ed{Th}

\Pf. ($\Rightarrow$) Since $(\mathcal{C}_{1},~\mathcal{C}_{2})$ and $(\mathcal{D}_{1},~\mathcal{D}_{2})$ are torsion pairs, we have that
$(\mathfrak{U}^{\mathcal{C}_{1}}_{\mathcal{D}_{1}})^{\perp_{0}}=\mathfrak{J}^{\mathcal{C}_{1}^{\perp_{0}}}_{\mathcal{D}_{1}^{\perp_{0}}}$
$=\mathfrak{J}^{\mathcal{C}_{2}}_{\mathcal{D}_{2}}$ by Proposition \ref{prop4} (2).
Also $^{\perp_{0}}(\mathfrak{J}^{\mathcal{C}_{2}}_{\mathcal{D}_{2}})=\mathfrak{U}^{^{\perp_{0}}\mathcal{C}_{2}}_{^{\perp_{0}}\mathcal{D}_{2}}=\mathfrak{U}^{\mathcal{C}_{1}}_{\mathcal{D}_{1}}$
by the proposition \ref{prop4} (1). So we have that $(\mathfrak{U}^{\mathcal{C}_{1}}_{\mathcal{D}_{1}},~\mathfrak{J}^{\mathcal{C}_{2}}_{\mathcal{D}_{2}})$ is a torsion pair.

($\Leftarrow$) (1) For any $X\in\mathcal{C}_{1}$ and $Y\in\mathcal{C}_{2}$, then $\binom{X}{0}\in\mathfrak{U}^{\mathcal{C}_{1}}_{\mathcal{D}_{1}}$
and $\binom{Y}{0}\in\mathfrak{J}^{\mathcal{C}_{2}}_{\mathcal{D}_{2}}$.
Since $(\mathfrak{U}^{\mathcal{C}_{1}}_{\mathcal{D}_{1}},~\mathfrak{J}^{\mathcal{C}_{2}}_{\mathcal{D}_{2}})$ is a torsion pair,
we have that $\mathrm{Hom}_{R}(X,~Y)\cong \mathrm{Hom}_{T}(\binom{X}{0},~\binom{Y}{0})=0$ by Lemma \ref{lem4} (1).
i.e., $\mathrm{Hom}_{R}(\mathcal{C}_{1},~\mathcal{C}_{2})=0$.

For any $M\in\mathcal{C}_{1}^{\perp_{0}}$, then $\binom{M}{0}\in\mathfrak{J}^{\mathcal{C}_{1}^{\perp_{0}}}_{\mathcal{D}_{1}^{\perp_{0}}}=
(\mathfrak{U}^{\mathcal{C}_{1}}_{\mathcal{D}_{1}})^{\perp_{0}}=\mathfrak{J}^{\mathcal{C}_{2}}_{\mathcal{D}_{2}}$ by Proposition \ref{prop4} (2). i.e., $M\in\mathcal{C}_{2}$.

For any $N\in$ $^{\perp_{0}}\mathcal{C}_{2}$, then $\binom{N}{0}\in\mathfrak{U}^{^{\perp_{0}}\mathcal{C}_{2}}_{^{\perp_{0}}\mathcal{D}_{2}}=$
$^{\perp_{0}}(\mathfrak{J}^{\mathcal{C}_{2}}_{\mathcal{D}_{2}})=\mathfrak{U}^{\mathcal{C}_{1}}_{\mathcal{D}_{1}}$ by Proposition \ref{prop4} (1). i.e., $N\in\mathcal{C}_{1}$.
So $(\mathcal{C}_{1},~\mathcal{C}_{2})$ is a torsion pair.

(2) For any $I\in\mathcal{D}_{1}$ and $J\in\mathcal{D}_{2}$, then $\binom{0}{I}\in\mathfrak{U}^{\mathcal{C}_{1}}_{\mathcal{D}_{1}}$
and $\binom{\mathrm{Hom}_{S}(U,~J)}{J}\in\mathfrak{J}^{\mathcal{C}_{2}}_{\mathcal{D}_{2}}$. By Lemma \ref{lem4}, we have that
$\mathrm{Hom}_{S}(I,~J)\cong \mathrm{Hom}_{T}(\binom{0}{I},~\binom{\mathrm{Hom}_{S}(U,~J)}{J})=0$
since $(\mathfrak{U}^{\mathcal{C}_{1}}_{\mathcal{D}_{1}},~\mathfrak{J}^{\mathcal{C}_{2}}_{\mathcal{D}_{2}})$ is a torsion pair.
i.e., $\mathrm{Hom}_{S}(\mathcal{D}_{1},~\mathcal{D}_{2})=0$.

For any $P\in\mathcal{D}_{1}^{\perp_{0}}$, then
$\binom{\mathrm{Hom}_{S}(U,~P)}{P}\in\mathfrak{J}^{\mathcal{C}_{1}^{\perp_{0}}}_{\mathcal{D}_{1}^{\perp_{0}}}=
(\mathfrak{U}^{\mathcal{C}_{1}}_{\mathcal{D}_{1}})^{\perp_{0}}=\mathfrak{U}^{\mathcal{C}_{2}}_{\mathcal{D}_{2}}$ by Proposition \ref{prop4} (2). i.e., $P\in\mathcal{D}_{2}$.

For any $Q\in$ $^{\perp_{0}}\mathcal{D}_{2}$, then $\binom{0}{Q}\in\mathfrak{U}^{^{\perp_{0}}\mathcal{C}_{2}}_{^{\perp_{0}}\mathcal{D}_{2}}=$
$^{\perp_{0}}(\mathfrak{J}^{\mathcal{C}_{2}}_{\mathcal{D}_{2}})=\mathfrak{U}^{\mathcal{C}_{1}}_{\mathcal{D}_{1}}$ by Proposition \ref{prop3} (2). i.e., $Q\in\mathcal{D}_{1}$.
So $(\mathcal{D}_{1},~\mathcal{D}_{2})$ is a torsion pair.
\ \hfill $\Box$

\vskip 10pt
By Proposition \ref{prop1}, Theorem \ref{th}, Theorem \ref{th1}, Theorem \ref{th2} and \cite[Lemma~2.3]{AMV}, we can get the following corollary immediately.

\bg{Cor}\rm{}
Let $\binom{A}{FA}\bigoplus\binom{0}{B}$ be silting (or, partial silting) in $T$-$\mathrm{Mod}$, $_{R}U$ be flat and $U_{S}$ be projective. Then the following conclusions hold.

$(1)$ If $S^{+}\in B^{\perp_{0}}$, then $(\mathfrak{B}^{\mathrm{Gen}A}_{A^{\perp_{0}}},~\mathfrak{U}^{\mathrm{Gen}B}_{B^{\perp_{0}}})$ is a torsion pair.

$(2)$ If $R\in\mathrm{Gen}A$, then $(\mathfrak{U}^{\mathrm{Gen}A}_{A^{\perp_{0}}},~\mathfrak{J}^{\mathrm{Gen}B}_{B^{\perp_{0}}})$ is a torsion pair.
\ed{Cor}

By Theorem \ref{2}, we can get the following result.
%
%
%

\bg{Cor}\rm{}\label{}
Taking $T=\left(\begin{array}{cc}R & 0 \\R & R \\\end{array}\right)$, then the statements hold.

(1) If $R^{+}\in\mathcal{D}$, then $(\mathcal{C},~\mathcal{D})$ is torsion pair in $R$-$\mathrm{Mod}$ if and only if
$(\mathfrak{B}^{\mathcal{C}}_{\mathcal{C}},~\mathfrak{U}^{\mathcal{D}}_{\mathcal{D}})$ is a torsion pair in $T$-$\mathrm{Mod}$.

(2) If $R\in\mathcal{C}$, then $(\mathcal{C},~\mathcal{D})$ is torsion pair in $R$-$\mathrm{Mod}$ if and only if
$(\mathfrak{U}^{\mathcal{C}}_{\mathcal{C}},~\mathfrak{J}^{\mathcal{D}}_{\mathcal{D}})$ is a torsion pair in $T$-$\mathrm{Mod}$.
\ed{Cor}

\section{Declarations}

\textbf{Conflict of interest} This work does not have any conflicts of interest.

{\small

}

\end{document}